\documentclass[11pt, twoside]{article}
\usepackage{bcc-cambridge-class}

\bcctitle{Asymptotic enumeration via graph containers and entropy}
\bccshorttitle{Asymptotic enumeration via graph containers and entropy}

\bccname{Jinyoung Park}
\bccshortname{Jinyoung Park} 

\bccaddressa{Department of Mathematics \\ Courant Institute of Mathematical Sciences, NYU \\ 251 Mercer Street \\ New York, New York 10012, USA}
\bccemaila{jinyoungpark@nyu.edu}

\usepackage{color}

\usepackage{graphicx}

\newcommand{\dist}{{\textrm{dist}}}
\newcommand{\mis}{{\textrm{mis}}}
\newcommand{\Hom}{{\textrm{Hom}}}
\newcommand{\Lip}{{\textrm{Lip}}}

\newtheorem{observation}[example]{Observation}

\begin{document}

\makebcctitle


\begin{abstract}
The container methods are powerful tools to bound the number of independent sets of graphs and hypergraphs, and they have been extremely influential in the area of extremal and probabilistic combinatorics. We will focus on more specialized graph container methods due to Sapozhenko (1987) that deal with sets in expander graphs. Entropy, first introduced by Shannon (1948) in the area of information theory, is a measure of the expected amount of information contained in a random variable. Entropy has seen lots of fascinating applications in a wide range of enumeration problems. In this survey article, we will discuss recent developments that exploit a combination of the two methods on enumerating graph homomorphisms.
\end{abstract}


\section{Introduction}\label{sec:intro}

\subsection{Framework}\label{sec:framework} Most of the enumeration problems discussed in this article are within the following framework. Consider a collection $\mathcal X$ (actually, a sequence of collections $\mathcal X_n$) whose size we want to estimate. Finding the exact size of $\mathcal X$ is hard, and in particular, there are no appropriate recurrences that we can leverage to apply a mathematical induction or a traditional generating function approach. Instead, there is a subcollection $\mathcal X' \subseteq \mathcal X$ that is easily understood. Our goal is to show that (upon finding the ``right" $\mathcal X'$) $\mathcal X'$ accounts for almost all $\mathcal X$, thus concluding that the size of $\mathcal X'$ in fact provides a good estimate for the size of $\mathcal X$. 

Of course, the tools to implement the above idea depend on the nature of the enumeration problem under consideration. In this article, we focus on problems related to the ``graph container method" and the ``entropy method." Before discussing these tools in more detail, we first introduce two classical enumeration problems to illustrate the abstract framework above: Dedekind's problem \cite{Dedekind1897Ueber} and calculating the number of independent sets in the hypercube \cite{Korshunov1983Number}, where the corresponding $\mathcal X$ are the collections of antichains in the Boolean lattice and independent sets in the hypercube, respectively.

\subsection{Dedekind's problem}\label{sec:Dedekind}

An \textit{antichain} of a poset $P=(X, \preceq)$ is a subset of $X$ in which any two elements are incomparable in the (partial) order $\preceq$. Dedekind's problem from 1897 asks for the total number $\alpha(n)$ of antichains in the $n$-dimensional Boolean lattice $\mathcal B_n=(2^{[n]}, \subseteq)$, where $[n]:=\{1, 2, \ldots, n\}$. (Equivalently, $\alpha(n)$ is the number of monotone Boolean functions $f:\{0,1\}^n \rightarrow \{0,1\}$, or the number of elements of the free distributive lattice on $n$ generators.) The Boolean lattice can be viewed as the poset on $\{0,1\}^n$ (the collection of binary strings of length $n$) with the partial order
\[x \preceq y ~\text{ if and only if } ~x_i \le y_i ~\text{ for all $i \in [n]$}\]
where $x_i$ denotes the $i$-th coordinate of $x$. In what follows, we will take the latter viewpoint.

It is easy to check that $\alpha(1)=3$ (the empty collection, the singletons $\{0\}$, and $\{1\}$) and $\alpha(2)=6$ (the empty collection, the singletons $\{00\}$, $\{01\}, \{10\}, \{11\}$, and $\{01, 10\}$), but there is no closed formula known for $\alpha(n)$ for general $n$. The sequence $\alpha(n)$ grows very fast and is not easy to compute even for small values of $n$. Indeed, even with modern computing power, the exact value of $\alpha(n)$ is known only for $n \le 9$, where the case of $n=9$ was obtained very recently by two independent groups of researchers \cite{Jakel2023Computation, VanHirtum2023Computation} more than 30 years after the determination of $\alpha(8)$ \cite{Wiedemann1991Computation}.

Having the discussion from Section~\ref{sec:framework} in mind, observe that there is a trivial lower bound on $\alpha(n)$:
\begin{equation}\label{eq:alpha_triv}\alpha(n) \ge 2^{{n \choose \lfloor n/2 \rfloor}}.\end{equation}
To see this, denote by $\mathcal L_i$ (the ``$i$-th layer" of $\mathcal B_n$) the set of elements (binary strings) in $\mathcal B_n$ that have precisely $i$ coordinates equal to 1, observing that each $\mathcal L_i$ ($i \in [0,n]$) is an antichain. In particular, by taking the subsets of $\mathcal L_{\lfloor n/2 \rfloor}$ (which are all antichains), we obtain precisely $2^{n \choose \lfloor n/2 \rfloor}$ many antichains.

Sperner's theorem \cite{Sperner1928Ein} -- which is usually considered as the first result in extremal set theory -- says that a largest antichain in $\mathcal B_n$ has size ${n \choose \lfloor n/2 \rfloor}$, and the only antichains to achieve this bound are $\mathcal L_{\lfloor n/2 \rfloor}$ and $\mathcal L_{\lceil n/2 \rceil}$ (the ``middle layer(s)" of $\mathcal B_n$). So the lower bound in \eqref{eq:alpha_triv} is the best one can obtain following the rather naive strategy of taking all subsets of some largest antichain.

Interestingly, this trivial lower bound turned out to be not too far from the truth. A breakthrough result of Kleitman \cite{Kleitman1969Dedekind} in 1969 pinpointed the asymptotics of the logarithm of $\alpha(n)$, showing that the exponent in the lower bound in \eqref{eq:alpha_triv} is actually asymptotically optimal for $\log\alpha(n)$ (in this article, $\log=\log_2$):

\begin{theorem}[Kleitman \cite{Kleitman1969Dedekind}] \label{thm:Kleitman}
    \[\log\alpha(n)<\left(1+O\left(\frac{\log n}{\sqrt n}\right)\right){n \choose \lfloor n/2 \rfloor}.\]
\end{theorem}
The $O\left(\frac{\log n}{\sqrt n}\right)$ term on the right-hand side was subsequently improved to $O\left(\frac{\log n}{n}\right)$ by Kleitman and Markowsky \cite{Kleitman1975Dedekinds}. The problem of finding asymptotics for $\log \alpha(n)$ has been revisited several times: the Kleitman-Markowsky bound was reproved by Kahn \cite{Kahn2001EntropyIndependent} using the entropy method, and by Balogh, Treglown, and Wagner \cite{Balogh2016Applications} using the container method. Pippenger \cite{Pippenger1999Entropy} obtained a similar bound with a slightly weaker error term based on the entropy method (whose detailed implementation was different from Kahn's approach).

In 1981, Korshunov \cite{Korshunov1981Number} determined $\alpha(n)$ itself up to a $(1+o(1))$ factor, thereby ``asymptotically solving" Dedekind's problem. Korshunov's paper was more than 100 pages, and a simpler -- though still involved -- proof was later given by Sapozhenko \cite{Sapozhenko1991Number}. In the statement below, we use $M={n \choose \lfloor n/2 \rfloor}$ for simplicity.

\begin{theorem}[Korshunov \cite{Korshunov1981Number}, Sapozhenko \cite{Sapozhenko1991Number}]
\label{thm:Dedekind}
    As $n\rightarrow \infty$, for $n$ even we have
    \[ \alpha(n)=(1+o(1))2^{M}\exp\left[\binom{n}{n/2+1}\left( 2^{-n/2}+\frac{n^2-2n-16}{32}\cdot2^{-n}\right)  \right],  \]
    while for $n$ odd
     \begin{center}
     $\displaystyle\alpha(n)=(2+o(1))2^{M}\cdot $ 
     \\\vspace{0.7 em}
     {\small $ \exp\left[\displaystyle\binom{n}{\frac{n+1}{2}} \left(2^{-(n+1)/2} + \frac{3n^2-19}{32}\cdot 2^{-(n+1)}\right)+\binom{n}{\frac{n+3}{2}}\left(\frac{1}{2}\cdot 2^{-(n+1)/2} + \frac{n+5}{8}\cdot 2^{-(n+1)}\right)\right].$}
     \end{center}
\end{theorem}

We remark that, while the above formulas might look technical, each of them is actually a natural asymptotic lower bound (for the corresponding parity) on $\alpha(n)$. That is, with writing $\mathcal A(\mathcal B_n)$ for the collection of antichains of $\mathcal B_n$, for each parity of $n$, there is a (simpler) subfamily $\mathcal A'_n \subseteq \mathcal A(\mathcal B_n)$ with $|\mathcal A'_n|$ asymptotically equal to the formula in Theorem~\ref{thm:Dedekind}. At the end of Section~\ref{sec:intro}, we will provide a description of $\mathcal A'_n$.

So how one proves Theorem~\ref{thm:Dedekind} is as explained in the framework in Section~\ref{sec:framework}. Korshunov and Sapozhenko showed that $|\mathcal A(\mathcal B_n) \setminus \mathcal A'_n|=o(|\mathcal A'_n|)$, thus concluding that $\mathcal A'_n$ accounts for almost all of the antichains in $\mathcal B_n$. A main tool in \cite{Sapozhenko1991Number} is the graph container method. Note that this proof also provides the structure of ``typical" antichains of $\mathcal B_n$, as the argument says that almost all antichains in $\mathcal B_n$ are the ones in the description of $\mathcal A'_n$.

\subsection{Independent sets in the hypercube}\label{sec:iQd}

The $n$-dimensional (discrete) \textit{hypercube} (or the \textit{Hamming cube}) $Q_n$ is the graph on the vertex set $\{0,1\}^n$ where two vertices are adjacent if and only if they differ in exactly one coordinate. Equivalently, $Q_n$ is the covering graph of the poset $\mathcal B_n$. Note that $Q_n$ is an $n$-regular bipartite graph with the unique bipartition classes $\mathcal E$ and $\mathcal O$, where $\mathcal E$ is the set of vertices (binary strings of length $n$) which have an even number of 1's in their coordinates. For simplicity, we will write $N$ for $2^{n}=|V(Q_n)|$. (So $|\mathcal E|=|\mathcal O|=N/2$.)

An \textit{independent set} in a graph is a set of vertices in which no two vertices are adjacent. Denote by $\mathcal I(G)$ the collection of independent sets in a graph $G$, and $i(G):=|\mathcal I(G)|$. We are interested in estimating $i(Q_n)$. This number again grows very fast, and a closed formula is not known.

As in Dedekind's problem, one can easily think of a trivial lower bound on this number:
\[i(Q_n) \ge 2\cdot2^{N/2}-1\]
because each subset of $\mathcal E$ or $\mathcal O$ is an independent set (and the empty set is counted twice). Korshunov and Sapozhenko \cite{Korshunov1983Number}, using the graph container method, showed that this trivial lower bound is not far from the truth. 

\begin{theorem}[Korshunov, Sapozhenko \cite{Korshunov1983Number}]\label{thm:iQd} As $n \rightarrow \infty$,
\begin{equation}\label{eq:iQd} i(Q_n)= 2(1+o(1))\sqrt e 2^{N/2}.\end{equation}    
\end{theorem}

As in Theorem~\ref{thm:Dedekind}, the right-hand side of \eqref{eq:iQd} is an easy asymptotic lower bound on $i(Q_n)$. In the next section, we discuss the construction of this lower bound. 

\subsubsection{Asymptotic lower bound on $i(Q_n)$}\label{Subsubsec:lb_iQd} 

Given an independent set $I$, we say a vertex $v$ is \textit{occupied} if $v \in I$. Let's informally say that an independent set of $Q_n$ is in a ``ground state" (or ``pure") if it is entirely contained in either $\mathcal E$ or $\mathcal O$. So there are two ground states ($\mathcal E$ and $\mathcal O$), and the number of independent sets in each ground state is precisely $2^{N/2}$. 

Next, we estimate the number of independent sets with some ``flaws" (or ``noises") and examine their contribution to $i(Q_n)$. Let's first consider independent sets $I$ that have exactly one vertex (a flaw) from $\mathcal O$.  The number of such independent sets is
\[N/2\cdot 2^{N/2-n}=2^{n-1}2^{N/2-n}=2^{-1}2^{N/2}\]
($N/2$ is the number of ways to choose a flaw $v \in \mathcal O$; once we pick $v$, then the $n$ vertices in $N(v) ~(\subseteq \mathcal E)$ must be unoccupied so we are left with $2^{N/2-n}$ choices for $I \cap \mathcal E$). In general, the number of independent sets $I$ that have exactly $k$ ``non-nearby" flaws from $\mathcal O$ (more precisely: we say $v$ and $w$ are \textit{nearby} if their neighborhoods intersect) is approximately
\begin{equation}\label{eq:k_flaws}{N/2 \choose k} \cdot 2^{N/2-kn}\sim \frac{(N/2)^k}{k!}2^{N/2-kn}=2^{N/2} \cdot \frac{2^{-k}}{k!}\end{equation}
as long as $k$ is ``not so large" ($N/2 \choose k$ approximates the number of ways to choose $k$ non-nearby flaws from $\mathcal O$; those $k$ vertices altogether have a neighborhood of size $kn$, so we have $2^{N/2-kn}$ choices left for $I \cap \mathcal E$).

Now, by summing \eqref{eq:k_flaws} over $k$ from $0$ to any (slow) $\omega(1)$ and multiplying it by 2 (since there are two ground states) we obtain the right-hand side of \eqref{eq:iQd} as an asymptotic lower bound on $i(Q_n)$.

\begin{rem}
    In connection with the framework in Section~\ref{sec:framework}, in the problem of estimating $i(Q_n)$, $\mathcal X'_n$ consists of the independents sets described in Section~\ref{Subsubsec:lb_iQd}, and Theorem~\ref{thm:iQd} asserts that $\mathcal X'_n$ accounts for almost all of $\mathcal I(Q_n)$. We reiterate that Theorem~\ref{thm:iQd} therefore also provides a \textit{structural} result -- it says that almost all independent sets in $Q_n$ have the structure described in Section~\ref{Subsubsec:lb_iQd}.
\end{rem}

\subsubsection{Towards the upper bound}\label{sec:iQd_ub} 

We first briefly motivate why it is ``sound" to expect that the independent sets in Section~\ref{Subsubsec:lb_iQd} should take up almost all of $\mathcal I(Q_n)$. In this section, we call an independent set \textit{bad} if $I$ is not of the type described in Section~\ref{Subsubsec:lb_iQd}. As a simplest example of bad independent sets, consider independent sets $I$ that have exactly two nearby flaws from $\mathcal O$ -- i.e., two occupied vertices in $\mathcal O$ whose neighborhoods overlap. The number of independent sets of this type is roughly (ignoring some irrelevant constant factors)
\begin{equation}\label{eq:nearby_flaws} N/2\cdot n^2\cdot 2^{N/2-2n}\approx n^2\cdot 2^{N/2-n}~(\ll 2^{N/2});\end{equation}
here $N/2$ is the number of ways to choose the first flaw $v$; now a crucial point is that once we pick $v$, the second flaw is within the second neighborhood of $v$ (since the two flaws are nearby), so there are at most $n^2$ choices for it;  the two flaws still have the total neighborhood of size roughly $2n$ (because a pair of vertices in $Q_n$ has at most two common neighbors), so we have roughly $2^{N/2-2n}$ choices for $I \cap \mathcal E$. 

So the contribution from \eqref{eq:nearby_flaws} is negligible (in fact, exponentially small) and absorbed in the $o(1)$ term in \eqref{eq:iQd}. Note that the reason why the overall contributions of the independent sets with nearby flaws and those of independent sets with non-nearby flaws are different is because the number of choices for two nearby flaws is significantly smaller than that for two non-nearby flaws, while the sizes of the neighborhoods of the flaws in the two cases remain more or less the same. 

While independent sets with a small number of flaws are not so hard to analyze (as we saw in the above computation), the story is completely different when we consider independent sets with a large set of flaws. Notice that the two competing factors in bounding the number of bad independent sets are 
\[\text{the number of choices for flaws vs. the size of the total neighborhood of the flaws.}\]
Informally, the former is the ``cost" for specifying flaws, and the latter is the ``penalty" from the (given) flaws. The argument surrounding \eqref{eq:nearby_flaws} illustrates the idea that in order to show the bad independent sets are negligible, we have to show that the penalty from the flaws of the independent sets exceeds their costs. So we need an effective upper bound on the number of subsets of $\mathcal O$ (or $\mathcal E$) that consist of ``nearby vertices" with a prescribed size of their neighborhood. The graph container lemma (Lemma~\ref{lem:graph_containers}) due to Sapozhenko \cite{Sapozhenko1987Number} serves this purpose, and this is one of the main focuses of this article.

The rest of the article is organized as follows. Section~\ref{sec:hom} defines the notion of graph homomophism and introduces some well-known examples. In Section~\ref{section:entropy}, we introduce some basics of the entropy function and provide an example to illustrate how the entropy method is used in solving enumeration problems. Section~\ref{sec:containers} starts from an overview of the container method, and then dives into a more specialized graph container lemma for expanders that can be used to get sharper results when it applies. Section~\ref{sec:asymp.hom} discusses sharp asymptotics for the number of graph homomorphisms, and hints at why a combination of the graph container method and the entropy method can be useful here. In Section~\ref{sec:t^n}, we discuss the problem of enumerating antichains of $[t]^n$, a generalized Boolean lattice, in connection with a first successful application of the graph container lemma for expanders in the ``irregular" setting. Section~\ref{sec:mis} discusses counting maximal independent sets in the hypercube, which leads to a development of an interesting isoperimetric inequality on the hypercube.

\begin{notn} Given a graph $G$, $E(G)$ and $V(G)$ denote the set of edges and vertices of $G$, respectively. For $u, v \in V(G)$, write $u \sim_G v$ (or $u \sim v$ if the host graph is clear) if $\{u, v\} \in E(G)$. For $U \subseteq V(G)$, $G[U]$ denotes the subgraph of $G$ induced by $U$. As usual, $N(x)$ denotes the neighborhood of $x$, and $N(U):=\cup_{x \in U} N(x)$. The \textit{codegree} of $x$ and $y$ ($x \ne y$) is $|N(x) \cap N(y)|$, and the codegree of $G$ is the maximum of the codgree of all pairs in $V(G)$. For $A,B \subseteq V(G)$, $\nabla(A,B):=\{\{a,b\} \in E(G):a \in A, b \in B\}$ and $e(A,B):=|\nabla(A,B)|$. $\dist_G(u,v)$ (or simply $\dist(u,v)$ if the host graph is clear) denotes the length of a shortest path in $G$ between $u$ and $v$. The \defword{diameter} of a graph $G$, denoted by $\text{diam}(G)$, is $\max_{u,v \in V(G)} \dist_G(u,v)$. 

For an integer $k$, we write $G^k$ for the $k$-th power of $G$, that is, $V(G^k)=V(G)$ and $u, v \in V(G^k)$ are adjacent iff $u \ne v$ and $\dist_G(u,v) \le k$.   We say a vertex set $A \subseteq V(G)$ is \defword{$k$-linked} if $G^k[A]$ is connected. The \textit{closure} of $A$, denoted by $[A]$, is $[A]:=\{v\in V(G):N(v)\subseteq N(A)\}$. 

For integers $a$ and $b$, $[a]:=\{1, 2, \ldots, a\}$ and $[a,b]:=\{a, a+1, \ldots, b\}$. We use ${[n] \choose k}$ to denote the family of all $k$-element subsets of $[n]$. 

We use standard asymptotic notation, and w.h.p. means ``with high probability." Logarithms are always base 2. 
\end{notn}

We close this section by giving the promised construction of the asymptotic lower bound for Dedekind's problem: we describe the family $\mathcal A'_n$ whose size is asymptotically equal to the (asymptotic) values of $\alpha(n)$ in Theorem~\ref{thm:Dedekind}. We only describe the construction for even $n$, for simplicity (the construction for odd $n$ is similar). Recall that $\mathcal L_i$ denotes the $i$-th layer of the Boolean lattice $\mathcal B_n$. We define the \textit{comparability graph} of $\mathcal B_n$ to be the graph on $V(\mathcal B_n)$ where two vertices are adjacent iff they are comparable. Observe that for any $k \in \mathbb Z^+$, a vertex set $A \subseteq V(G)$ can be decomposed into its maximal $k$-linked subsets, and we call these subsets the \textit{$k$-linked components} of $A$ (with respect to $G$).

    \begin{theorem}\label{thm:Dedekind_structure} Let $n$ be even. Let $\mathcal A'_n$ be the collection of the antichains $I$ of $\mathcal B_n$ that satisfy the following properties:
\begin{enumerate}[(i)]
    \item $I \subseteq \bigcup_{n-1 \le i \le n+1} \mathcal L_i$;
    \item every 2-linked component of $I \cap (\mathcal L_{n-1} \cup \mathcal L_{n+1})$ with respect to the comparability graph of $\mathcal B_n$ is of size at most 2.
\end{enumerate}
Then, as $n \rightarrow \infty$,
\[\alpha(n)=\left(1+e^{-\Omega(n)}\right)|\mathcal A_n'|.\]    \end{theorem}
We leave it to the readers that $|\mathcal A_n'|$ (for even $n$) is asymptotically equal to the value of $\alpha(n)$ in Theorem~\ref{thm:Dedekind}.

\section{Graph homomorphisms}\label{sec:hom}

Independent sets and many other important combinatorial objects can be incorporated into the notion of graph homomorphisms which we introduce now. Given graphs $G$ and $H$, a \textit{graph homomorphism} from $G$ to $H$ (or an \textit{$H$-coloring} of $G$) is a map $f:V(G) \rightarrow V(H)$ such that if $\{v, w\} \in E(G)$, then $\{f(v), f(w)\} \in E(H)$. Write $\Hom(G,H)$ for the collection of graph homomorphisms from $G$ to $H$. Below are some well-known graph homomorphisms.

\begin{example}[Independent set] Let $H$ be the graph with $V(H)=\{v_0, v_1\}$ and $E(H)=\{\{v_0, v_0\}, \{v_0, v_1\}\}$.
\begin{figure}[h]
\begin{center}
\includegraphics[height=1cm]{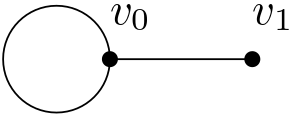}
\end{center}
\caption{$H$ for independent sets}
\end{figure}

To see the correspondence between an independent set and an $H$-coloring, for a given independent set $I \in \mathcal I(G)$, define $f_I$ to be
\[f_I(v)=\begin{cases}v_0 & \text{ if $v \not\in I$;}\\
v_1 & \text{ if $v \in I$}.\end{cases}\]
Then $f_I \in \Hom(G,H)$, and it is easy to see that the above correspondence between $I$ and $f_I$ provides a bijection between $\mathcal I(G)$ and $\Hom(G,H)$. That is, any independent set can be regarded as an $H$-coloring, and vice versa.

\end{example}

\begin{example}[Proper coloring] Set $H=K_q$, the clique on $q$ vertices. It is easy to see that there is a one-to-one correspondence between $H$-colorings of $G$ and (proper) $q$-colorings of $V(G)$. 
\end{example}

We will drop ``proper" and just say a $q$-coloring or a coloring for a proper coloring. In the rest of the paper, we will use $\mathcal C_q(G)$ for the collection of $q$-colorings of $V(G)$, and $c_q(G):=|\mathcal C_q(G)|.$

\begin{example}[Lipschitz function] Given an integer $M$, say an integer-valued function $g: V(G) \rightarrow \mathbb Z$ is \textit{$M$-Lipschitz} if $v \sim w$ implies $|g(v)-g(w)|\le M$. 

Let $H$ be the graph with $V(H)=\mathbb Z$ and $E(H)=\{\{m,m\}, \{m, m+1\}:m \in \mathbb Z\}$ as in the figure below.
\begin{figure}[h]
\begin{center}
\includegraphics[height=1.2cm]{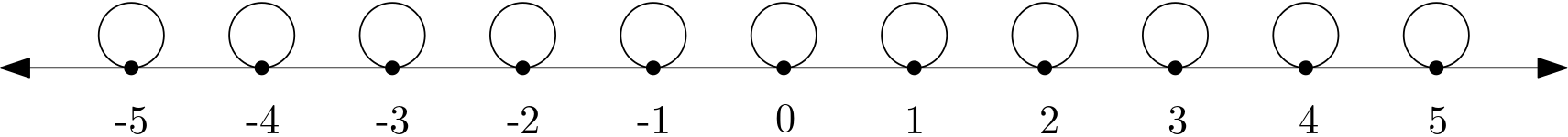}
\end{center}
\caption{$H$ for 1-Lipschitz functions}
\end{figure}

\noindent Observe that there is a one-to-one correspondence between $H$-colorings of $G$ and 1-Lipschitz functions from $G$. Similarly, if we define $H_M$ to be the graph on $\mathbb Z$ with $E(H_M)=\{\{m,m\}, \{m, m+1\}, \ldots, \{m, m+M\}:m \in \mathbb Z\}$, then there is a bijection between $H_M$-colorings of $G$ and $M$-Lipschitz functions from $G$.
\end{example}

As a last example of this section, an integer-valued function $g:V(G) \rightarrow \mathbb Z$ is called a \textit{$\mathbb Z$-homomorphism} if $v \sim w$ implies $|g(v)-g(w)|=1$. A $\mathbb Z$-homomorphism is an $H$-coloring where $H$ is the integer line. Note that $G$ must be bipartite to admit a $\mathbb Z$-homomorphism.

\section{Entropy}\label{section:entropy}

The entropy of a random variable $\mathbf X$ was first introduced by Shannon in 1948 \cite{Shannon1948Mathematical} as a measure of the ``information" or ``uncertainty" contained in $\mathbf X$. Entropy has been used to tackle many enumeration problems, and this article focuses on its usefulness for counting (or bounding the number of) independent sets, colorings, or more generally, graph homomorphisms. In this section, we quickly introduce some very basics of the entropy function, and provide an example to show how entropy is used for enumerating proper colorings of a given graph.

\subsection{Basics}\label{sec:entropy} This section briefly collects some basics of entropy; for a less hurried introduction, see, e.g., McEliece \cite{McEliece1977Theory}.

For a discrete {random variable} $\mathbf{X}$, the \textit{entropy} of $\mathbf X$ is
    $$H(\mathbf{X})=\sum_x \mathbb P(\mathbf X=x)\log \frac{1}{\mathbb P(\mathbf X=x)}$$
    where the sum is taken over all $x$ with $\mathbb P(\mathbf X =x) \ne 0$. (By convention, we take $0 \log \frac{1}{0}=0$.) The entropy of $\mathbf X$ is often understood as the amount of ``randomness" of $\mathbf X$; for example, one can show that
    \begin{equation}\label{eq:entropy_a} H(\mathbf{X})\leq \log |\text{Image}(\mathbf{X})|\end{equation}
    (Image($\mathbf X$) is the set of values that $\mathbf X$ takes with nonzero probability) by Jensen's inequality using the concavity of the logarithm function. Furthermore, in \eqref{eq:entropy_a} the equality holds iff $\mathbf{X}$ is chosen uniformly at random from Image($\mathbf{X}$). So roughly,  the ``more random" the random variable $\mathbf X$ is, the larger the entropy of $\mathbf X$ is.
    
For any event $T$, we define the entropy of $\mathbf{X}$ given $T$ as
\[H(\mathbf{X}|T)=\sum_x \mathbb P(\mathbf X=x|T)\log \frac{1}{\mathbb P(\mathbf X=x|T)}.\]
    For another discrete random variable $\mathbf Y$, the \textit{conditional entropy} of $\mathbf{X}$ with respect to  $\mathbf{Y}$ is 
    \begin{equation}\label{eq:cond_entropy} \begin{split}H(\mathbf{X}|\mathbf{Y})&=\sum_y \mathbb P(\mathbf Y=y) H(\mathbf{X}|\mathbf{Y}=y)\\
&= \sum_y \mathbb P(\mathbf Y=y) \sum_x \mathbb P(\mathbf X=x|\mathbf Y=y)\log \frac{1}{\mathbb P(\mathbf X=x|\mathbf Y=y)}.\end{split}\end{equation}

    Below are some standard facts about the entropy; $\mathbf X$, $\mathbf Y$, and $\mathbf X_i$'s are all used for discrete random variables. For any random vector $\mathbf{X}=(\mathbf{X}_1,\dots,\mathbf{X}_n)$ we have
    \begin{equation}\label{eq:entropy_d} H(\mathbf{X})\leq H(\mathbf{X}_1)+\dots+H(\mathbf{X}_n).\end{equation}
The above inequality can be sharpened using the conditional entropy; for example, for two random variables, we have
    \begin{equation}\label{eq:entropy_b} H(\mathbf{X},\mathbf{Y})=H(\mathbf{X})+H(\mathbf{Y}|\mathbf{X}).\end{equation}
It also holds that
    \begin{equation}\label{eq:entropy_c} {\text{if $\mathbf{Y}$ determines $\mathbf{X}$, then $H(\mathbf Z|\mathbf{Y}) \le H(\mathbf Z|\mathbf{X});$}}\end{equation}    
intuitively, the above means that the amount of information in $\mathbf Z$ that remains after revealing $\mathbf X$ is at least the amount remaining after revealing $\mathbf Y$.

The following result is a celebrated inequality due to Shearer, which generalizes the subadditivity property \eqref{eq:entropy_d}.

	{\begin{lemma}[\cite{Chung1986Some}]\label{lem:Sh}
		Let $\mathbf{X}=(\mathbf{X}_1, \ldots, \mathbf{X}_k)$ be a random vector, and let $\alpha_A\in\mathbb R^+$ be given for every $A\subseteq [k]$. If $\sum_{A \ni i} \alpha_A\ge 1$ for all $i \in [k]$,
		then
		\[H(\mathbf{X})\le \sum_{A \subseteq [k]} \alpha_A H(\mathbf{X}_A),\]
		where $\mathbf{X}_A=(\mathbf{X}_i:i \in A)$.
	\end{lemma}}

\subsection{Enumeration via entropy}  In this section, we provide a standard-ish entropy-based argument for bounding the number of colorings of regular bipartite graphs. The proof below, which is a special case of the results in \cite{Galvin2004Weighted}, can easily be extended to more general graph homomorphisms. Recall that $\mathcal C_q(G)$ is the collection of $q$-colorings of $V(G)$, and $c_q(G):=|\mathcal C_q(G)|.$

\begin{prop}\label{prop:entropy_coloring}
        Let $G$ be a $d$-regular bipartite graph on $m$ vertices. For any $q$,
        \[\log c_q(G)\le \frac{m}{2}\cdot \left(\log\left(\left\lfloor \frac{q}{2} \right\rfloor \left\lceil \frac{q}{2}\right\rceil\right)+O_q(1/d)\right)\]
where in $O_q(\cdot)$, the implicit constant (only) depends on the value of $q$.
\end{prop}

Observe that, trivially, $\log c_q(G) \ge \frac{m}{2}\cdot \log\left(\left\lfloor \frac{q}{2} \right\rfloor \left\lceil \frac{q}{2}\right\rceil\right)$ (fix $\lfloor \frac{q}{2} \rfloor$ colors and call this set of colors $\mathcal C$; denote the set of the remaining $\lceil \frac{q}{2} \rceil$ colors by $\mathcal C'$ and consider every coloring of $G$ (with a bipartition $(A,B)$) where each vertex in $A$ is given a color from $\mathcal C$ and each vertex in $B$ is given a color from $\mathcal C'$). So if $q=O(1)$, then Proposition~\ref{prop:entropy_coloring} provides an asymptotically tight upper bound on $\log c_q(G)$.

We prove the above statement only for $q=4$ for a simpler exposition, so our target bound is
\begin{equation}\label{eq:c4.target} \log c_4(G) \le (1+O(1/d))m.\end{equation}
 The proof below extends almost identically to any $q$.

\begin{proof} Denote by $X$ and $Y$ a bipartition of $G$, noting that $|X|=|Y|=m/2$. We write $N_x$ for $N(x)$. For $f \in \mathcal C_4(G)$ and $x \in V(G)$, $f_x$ denotes the value of $f$ at $x$, and for $A \subseteq V(G)$, $f_A=(f_x:x \in A)$ and $f(A)=\{f_x:x \in A\}$ (that is, $f_A$ contains the full information about the color assignment of each $x \in A$, while $f(A)$ (the ``palette" of $A$) only contains the information about what colors are used in $A$).

Let $\mathbf f$ be an element of $\mathcal C_4(G)$ chosen uniformly at random. A key point is that, by the equality condition of \eqref{eq:entropy_a},
\[H(\mathbf f)=\log c_4(G).\]
Therefore, the conclusion \eqref{eq:c4.target} will follow if we show that $H(\mathbf f)$ is bounded above by the right-hand side of \eqref{eq:c4.target}. To this end, first, by \eqref{eq:entropy_b}, 
\begin{equation}\label{eq:entropy_coloring1} H(\mathbf f)=H(\mathbf f_Y)+H(\mathbf f_X|\mathbf f_Y).\end{equation}
The second term of the right-hand side of the above is, using a conditional version of \eqref{eq:entropy_d} and then \eqref{eq:entropy_c},
\[H(\mathbf f_X|\mathbf f_Y) \le \sum_{x \in X} H(\mathbf f_x|\mathbf f_Y)\le \sum_{x \in X}H(\mathbf f_x|\mathbf f(N_x));\]
the first term of the right-hand side of \eqref{eq:entropy_coloring1} is (a justification will follow)
\[H(\mathbf f_Y) \le \frac{1}{d}\sum_{x \in X}H(\mathbf f_{N_x})=\frac{1}{d}\sum_{x \in X}H(\mathbf f(N_x),\mathbf f_{N_x})=\frac{1}{d}\sum_{x \in X} \left[H(\mathbf f(N_x))+H(\mathbf f_{N_x}|\mathbf f(N_x))\right]\]
where for the inequality, we use Lemma~\ref{lem:Sh} with
    \[\alpha_A=\begin{cases} 1/d & \text{ if } A=N_x \text{ for some } x \in X \\ 0 & \text{ otherwise} \end{cases}\]
(and the second equality uses \eqref{eq:entropy_b}). Combining the two bounds, \eqref{eq:entropy_coloring1} is
\[\begin{split}H(\mathbf f)&\le \sum_{x \in X} \frac{1}{d}\left[H(\mathbf f(N_x))+H(\mathbf f_{N_x}|\mathbf f (N_x))\right]+H(\mathbf f_x|\mathbf f(N_x)).\end{split}\]
Let $T(x)$ ($x \in X$) be the summand of the above sum, that is,
    \[\begin{split}T(x)=\frac{1}{d}\left[H(\mathbf f(N_x))+H(\mathbf f_{N_x}|\mathbf f(N_x))\right]+H(\mathbf f_x|\mathbf f(N_x)).\end{split}\]
    We claim that $T(x)\le 2+O(1/d)$ from which the conclusion will follow. We will use (a conditional version of) \eqref{eq:entropy_a} repeatedly. The first term of the above expression is easily seen to be $O(1/d)$, since when there are only 4 colors, $H(\mathbf f(N_x))=O(1)$ (the number of possible palettes for $N_x$ cannot exceed $2^4$). To bound the last two terms, note that for each possible value (a set of colors) $C$ of $\mathbf f(N_x)$, by \eqref{eq:entropy_a},
    \[H(\mathbf f_x|\mathbf f (N_x)=C) \le \log(4-|C|); \text{ and}\]
    \[H(\mathbf f_{N_x}|\mathbf f(N_x)=C)\le \sum_{y \in N_x} H(\mathbf f_y|\mathbf f(N_x)=C) \le d\log|C|.\]
    Since $\log z+\log(4-z) \le 2$ for any $z$, this bounds $\frac{1}{d}H(\mathbf f_{N_x}|\mathbf f(N_x))+H(\mathbf f_x|\mathbf f(N_x))$ by (using the definition of conditional entropy \eqref{eq:cond_entropy})
    \[\sum_C \mathbb P(\mathbf f(N_x)=C)\left[\frac{1}{d}H(\mathbf f_{N_x}|\mathbf f(N_x)=C)+H(\mathbf f_x|\mathbf f(N_x)=C)\right] \le 2,\]
    yielding $T(x) \le 2+O(1/d)$ for any $x \in X$. \end{proof}

\section{The method of graph containers and sharp asymptotics}\label{sec:containers}

The container method has been proven to be an essential tool in the study of independent sets in graphs and hypergraphs.  In the graph setting, its roots can be traced back to the work of Kleitman and Winston~\cite{Kleitman1982Number} and Sapozhenko~\cite{Sapozhenko1987Number}. These tools were extended to the context of hypergraphs and developed into a powerful method in the seminal works of Balogh, Morris and Samotij~\cite{Balogh2015Independent} and Saxton and Thomason~\cite{Saxton2015Hypergraph}. The hypergraph container method has applications in broad areas including  extremal graph theory, Ramsey theory, additive combinatorics, and discrete geometry. We refer the readers to \cite{Balogh2018Method} for an overview of the hypergraph container method.

A main focus of this article is a more specialized graph container lemma due to Sapozhenko.  As a warm-up, in Section~\ref{sec:containers_warmup}, we introduce some basic concepts of the general graph container method. In Section~\ref{sec:graph_containers}, we provide the statement of Sapozhenko's graph container lemma for expanders. This is a powerful tool for the problems of obtaining asymptotically sharp enumeration results and thus gaining an understanding of the structure of typical members of the family under consideration.

\subsection{Warm-up}\label{sec:containers_warmup}

The method of graph containers provides an efficient upper bound on the number of certain vertex subsets (often independent sets) of a given graph.

\begin{Def}[System of containers] Given a graph $G$ and $\mathcal F \subseteq 2^{V(G)}$, suppose a family $\mathcal C \subseteq 2^{V(G)}$ satisfies the following property:
\[\text{for any $F \in \mathcal F$, there is $C \in \mathcal C$ such that $F \subseteq C$.}\]
An element of the family $\mathcal C$ is called a \textit{container}, and the family $\mathcal C$ is called a \textit{system of containers} for $\mathcal F$.     
\end{Def}

Let some family $\mathcal F$ be given, and suppose our goal is to approximate $|\mathcal F|$. In the approach using the container method, we are hoping that there is a system of containers $\mathcal C$ for $\mathcal F$ such that both $|\mathcal C|$ and $\max\{|C|:C \in \mathcal C\}$ are reasonably small (or close to the truth). Then the naive inequality
\begin{equation}\label{eq:containers} |\mathcal F| \le \sum_{C \in \mathcal C}2^{|C|} \le |\mathcal C|2^{\max\{|C|:C \in \mathcal C\}}\end{equation}
will produce a good bound.
This approach turned out to work incredibly well (modulo the noble ideas in the works referenced at the beginning of Section~\ref{sec:containers}) for the family of independent sets in graphs (and certain hypergraphs). A very rough intuition behind this approach is that often independent sets in a graph cluster together, so specifying independent sets as subsets of a container should provide an efficient counting.

Below is one of the very first applications of the graph container method due to Sapozhenko \cite{Sapozhenko2001Number}.

\begin{theorem}\label{thm:Sap_indep}
    Let $G$ be a $d$-regular graph on $m$ vertices. Then
    \[\log i(G)<\left(1+O\left(\sqrt\frac{\log d}{ d}\right)\right)\frac{m}{2}.\]
\end{theorem}

    Theorem~\ref{thm:Sap_indep} follows from the following lemma that guarantees the existence of a suitable system of containers. (See \cite[Theorem 1.6.1]{Alon2016Probabilistic} for its proof.)

    \begin{lemma}\label{lem:Sap_indep} Let $G$ be as in Theorem~\ref{thm:Sap_indep}. 
        For any $\epsilon>0$, there is a collection $\mathcal C \subseteq 2^{V(G)}$ with the following properties:

        (i) For each independent set $I$ in $G$, there is a member $C \in \mathcal C$ such that $I \subseteq C$;

        (ii) For every $C \in \mathcal C$, $|C| \le \frac{m}{\epsilon d}+\frac{m}{2-\epsilon}$;

        (iii) $|\mathcal C|\le \sum_{i \le m/(\epsilon d)} {m \choose i}$.
    \end{lemma}

To see that Theorem~\ref{thm:Sap_indep} follows from the above lemma, using \eqref{eq:containers} (with $\mathcal F=\mathcal I(G)$),
\[i(G) \le \sum_{i \le m/(\epsilon d)} {m \choose i} \cdot \exp_2\left[{\frac{m}{\epsilon d}+\frac{m}{2-\epsilon}}\right].\]
By the well-known estimate that $\log\left[\sum_{i \le m/(\epsilon d)}{m \choose i}\right] = O\left(\frac{m\log(\epsilon d)}{\epsilon d}\right),$ we have
\[\log i(G)=O\left(\frac{m \log(\epsilon d)}{\epsilon d} + \frac{m}{\epsilon d}+\frac{m}{2-\epsilon}\right).\]
The right-hand side of the above is minimized when $\epsilon=\Theta(\sqrt{\log d/d})$. By plugging in this value for the above expression, we obtain the conclusion of Theorem~\ref{thm:Sap_indep}.

\subsection{Graph container lemma for expanders}\label{sec:graph_containers}

In the previous section, we discussed an upper bound on $\log |\mathcal F|$ for $\mathcal F=\mathcal I(G)$ using the method of graph containers. Notice that the above argument is not enough to give asymptotics for $|\mathcal F|$ itself. On the other hand, Sapozhenko showed that, if given graphs satisfy certain regularities in degrees and exhibit certain expansion properties, then one can further refine the container argument and obtain sharp asymptotics (as in Theorems~\ref{thm:Dedekind} and \ref{thm:iQd}). In this section, we state a graph container lemma for expanders\footnote{This name can be misleading --  the expansion that the lemma requires is much smaller than the usual notion of expanders.} (Lemma~\ref{lem:graph_containers}) from \cite{Sapozhenko1987Number}, and sketch how to derive the sharp asymptotics in Theorem~\ref{thm:iQd} from this lemma. Recall the definition of 2-linkedness from the notation section in Section \ref{sec:intro}.

\begin{Def}
    Let $G$ be a $d$-regular, bipartite graph with bipartition classes $X$ and $Y$. For positive integers $a$ and $g$ and $v \in Y$, define
\[\mathcal G(a,g, v)=\{A \subseteq X: \text{$A$ is 2-linked}, |[A]|=a, |N(A)|=g, v \in N(A)\}.\]
(Recall that $[A]$ denotes the closure of $A$.)
\end{Def}

\begin{lemma}\label{lem:graph_containers}
    For each pair of constants $c_1, c_2>0$, there is a constant $c=c(c_1, c_2)>0$ such that the following holds. If $G$ is a $d$-regular bipartite graph with partition classes $X$ and $Y$ and codegree $c_1$, then for $g \ge d^4$ and $g-a>\frac{c_2 g \log^3 d}{d^2}$, we have
    \begin{equation}\label{eq:agv}|\mathcal G(a, g, v)| \le 2^{g-c(g-a)/\log d}.\end{equation}
\end{lemma}

Sapozhenko's proof \cite{Sapozhenko1987Number} for the above lemma uses a container-type argument (but different from that for Lemma~\ref{lem:Sap_indep}). The original proof  is written in Russian, and Galvin \cite{Galvin2019Independent} provides a more accessible exposition of the proof of Theorem~\ref{thm:iQd}, including a proof of Lemma~\ref{lem:graph_containers}.

\begin{rem}
1. Sapozhenko's original formulation of the above lemma doesn't require strict regularity (as long as certain expansion requirements similar to that in Lemma~\ref{lem:graph_containers} are satisfied). For example, the degree conditions required in \cite{Sapozhenko1987Number} are the following: there is a constant $C$ such that
\[\text{$\max_{v \in X} \deg(v)\le C\cdot \min_{v \in X} \deg(v)$; \quad $\max_{v \in Y} \deg(v)\le C\cdot \min_{v \in Y} \deg(v)$,}\]
and most crucially,
    \begin{equation}\label{eq:Sap_deg} \max_{v \in Y} \deg(v) \le \min_{v \in X} \deg(v).\end{equation}
A proof of the lemma under the above (more relaxed) degree conditions can be found  from   \cite{Jenssen2024Refined} in a more general hard-core model setting.

2.  The exponential term $-\Omega((g-a)/\log d)$ on the right-hand side of \eqref{eq:agv} is improved to $-\Omega(g-a)$ in \cite{Park2022Note} (and implicitly in \cite{Kahn2022Number}). The weaker bound $-\Omega((g-a)/\log d)$ was already enough for many applications, but the improved bound $-\Omega(g-a)$ was crucial for obtaining sharp estimates for more restricted classes of independent sets including maximal independent sets in $Q_n$ \cite{Kahn2022Number} and ``balanced independent sets" in $Q_n$ \cite{Park2022Note}.    

3. Jenssen and Perkins \cite{Jenssen2020Independent} showed that a combination of the graph container lemma for expanders (more precisely, its variant for the hard-core model by Galvin \cite{Galvin2011Threshold}) and the ``cluster expansion method" from statistical physics results in powerful consequences. One consequence is the following vast refinement of the asymptotics in Theorem~\ref{thm:iQd}:
\[2\sqrt e\cdot 2^{N/2}\left(1+\frac{3n^2-3n-2}{8 N}+\frac{243n^4-646n^3-33n^2+436n+76}{384 N^2}+O\left(\frac{n^6}{N^{3}}\right)\right).\]
In fact, Jenssen and Perkins provided a formula and an algorithm for computing the asymptotics of $i(Q_n)$ to arbitrary order in $N^{-1}$.
We refer the readers to \cite{Jenssen2020Independent} for the further fruitful applications of this combination of methods. For the cluster expansion method itself, we refer to a survey paper by Jenssen \cite{Jenssen2024Cluster}. Jenssen and Perkins' idea of combining the two methods has been further adapted for lots of enumeration problems, including \cite{Balogh2021Independent, Collares2025Counting, Geisler2025Counting, Jenssen2023Homomorphisms, Jenssen2024Dedekinds, Jenssen2025Number, Kronenberg2022Independent, Li2025Number}. Jenssen, Malekshahian, and the author \cite{Jenssen2024Dedekinds} used the cluster expansion method to obtain a similar refinement of the asymptotics in Theorem~\ref{thm:Dedekind}. For example, for even $n$,
\[\alpha(n)=\left(1+o\left(2^{-2n}n^{O(1)}\right)\right) {2^{\binom{n}{n/2}}} \times\]
{\small
    \[{\exp\left[\binom{n}{\frac{n}{2}-1}\left(2^{-\frac{n}{2}}+(n^2-2n-16)2^{-n-5}+\frac{1}{3}(3n^4-4n^3-60n^2+112n+512)2^{-\frac{3n}{2}-9}\right) \right]}.\]}

\noindent Again, as in \cite{Jenssen2020Independent}, we also provide a formula and an algorithm for computing the asymptotics of $\alpha(n)$ to arbitrary order in $2^{-n/2}$.
\end{rem}
    
Now we sketch the derivation of Theorem~\ref{thm:iQd} from Lemma~\ref{lem:graph_containers}. For the sake of clarity, in the example below, we consider only a special type of bad independent set (recall from Sections \ref{Subsubsec:lb_iQd} and \ref{sec:iQd_ub} that bad independent set means it contains at least 2 nearby flaws) in $\mathcal I(Q_n)$ and show that their contribution to $i(Q_n)$ is negligible using Lemma~\ref{lem:graph_containers}. The (full) derivation of Theorem~\ref{thm:iQd} essentially follows by summing up the contributions of all types of bad independent sets. The contribution of each type of bad independent set is exponentially small, so the total contribution of bad independent sets to $i(Q_n)$ is again negligible. A detailed computation can be found in \cite[Corollary 4.1]{Galvin2019Independent}.

\begin{example}\label{ex:iQd}
Recall that $Q_n$ has the bipartition class $\mathcal E$ and $\mathcal O$ and $N=2^n=|V(Q_n)|$ (so $|\mathcal E|=|\mathcal O|=N/2$).

Let $\mathcal I_k \subseteq \mathcal I(Q_n)$ be the collection of independent sets $I$ such that $I \cap \mathcal O$ is 2-linked and $|[I \cap \mathcal O]|=k$ where $n^4 \le k \le N/4$. (We pose the lower bound on $k$ so that Lemma~\ref{lem:graph_containers} is applicable; the upper bound $N/4$ is natural since if $I \in \mathcal I(Q_n)$, then $\min\{|[I \cap \mathcal E]|, |[I \cap \mathcal O]|\} \le N/4$ by the expansion property of $Q_n$ that for any $A \subseteq \mathcal E$ or $A \subseteq \mathcal O$, $|A| \le |N(A)|$.) We claim that $|\mathcal I_k|$ is exponentially smaller than $2^{N/2}$ (the main term of the (good) independent set count in Section~\ref{Subsubsec:lb_iQd}). To see this, note that
\[\begin{split}|\mathcal I_k|&= \sum_{\substack{A \subseteq \mathcal O\\ A \text{ 2-linked}, |[A]|=k}} 2^{N/2-|N(A)|}\\
 &\le 2^{N/2}\sum_{v \in \mathcal E}\sum_{g} \sum_{\substack{A \subseteq \mathcal O, v \in N(A)\\ A \text{ 2-linked}, |[A]|=k, |N(A)|=g}}2^{-|N(A)|}\\
& = 2^{N/2}\sum_{v \in \mathcal E} \sum_{g} |\mathcal G(k, g, v)|\cdot 2^{-g}.\end{split}\]
A well-known isoperimetric inequality on the hypercube (see, e.g., \cite[Claim 2.5]{Galvin2019Independent}) says that if $Y \subseteq \mathcal O$ satisfies $|Y|\le N/4$, then $|N(Y)|-|Y|=\Omega(|N(Y)|/\sqrt n)$. Therefore, by applying Lemma~\ref{lem:graph_containers} (and applying the isoperimetric inequality with $Y=[A]$), 
$|\mathcal G(k,g,v)|\cdot 2^{-g}$ is at most $2^{-\Omega(g-k)/\log n}=2^{-\Omega(g/(\sqrt n\log n))}$ if $g$ satisfies $g-k=\Omega(g/\sqrt n)$; and 0 otherwise. Now, the last line of the above display is at most
\[2^{N/2}\cdot N/2 \cdot \sum_{g \ge k} 2^{-\Omega(g/(\sqrt n\log n))}=2^{N/2} \cdot 2^{-\Omega(k/(\sqrt n\log n))}\]
since $k \ge n^4$. This yields the desired conclusion.
\end{example}

\subsection{Applications to independent set counts} Lemma~\ref{lem:graph_containers} and its variants have been used for various host graphs to obtain sharp asymptotics for the number of independent sets and/or the structure of typical independent sets. The studied graphs include the discrete torus $\mathbb Z^n_m$ for even $m$ by Jenssen and Keevash \cite{Jenssen2023Homomorphisms}, the middle layers of the hypercube by Balogh, Garcia, and Li \cite{Balogh2021Independent}, and Abelian Cayley graphs by Potukuchi and Yepremyan \cite{Potukuchi2021Enumerating}. Recently, Collares, Erde, Geisler, and Kang \cite{Collares2025Counting} extended this approach to a broader class of regular bipartite graphs including the Cartesian product of bipartite, regular base graphs of bounded size. Using the observation that intersecting families are independent sets in a Kneser graph, Balogh, Garcia, Li, and Wagner \cite{Balogh2021Intersecting} showed that almost all $k$-uniform intersecting families in $2^{[n]}$ are ``stars" if  $n \ge 2k+100\ln k$. Their proof uses the graph container method and the Das--Tran removal lemma \cite{Das2016Removal}.

The entropy method, graph container method, and/or their combinations also apply to the setting of the ``hard-core model," a fundamental model in statistical physics that generalizes the notion of independent sets. The present article will not survey the huge literature on this topic, and just refer to \cite{Kahn2001Entropy, Galvin2011Threshold, Jenssen2024Refined} as indicators to the works more closely related to the current discussions.

The independent set enumeration problems have also been considered in sparse random settings. For example, Kronenberg and Spinka \cite{Kronenberg2022Independent} obtained asymptotics for the number of independent sets in edge-percolated hypercubes by studying an antiferromagnetic Ising model on the hypercube. Geisler, Kang, Sarantis, and Wdowinski \cite{Geisler2025Counting} extended Kronenberg-Spinka's approach to a broader class of graphs, including even tori of growing side-length. Several cornerstone results in extremal combinatorics such as Erd\H{o}s-Ko-Rado theorem \cite{Erdos1961Intersection} and Sperner's theorem \cite{Sperner1928Ein} can also be stated using the language of independent sets. Sparse random versions of those problems were studied in \cite{Balogh2014Random, Collares2016Maximum, Hamm2019Erdos, Balogh2023Sharp}, where the approaches based on the graph container methods were again a key ingredient.

\section{Sharp asymptotics for the number of graph homomorphisms}\label{sec:asymp.hom}

In this section, our main concern is getting sharp asymptotics for the number of graph homomorphisms from given graphs. Again, this will also provide us with the typical structure of the graph homomorphisms under consideration. 

We start our discussion with graph colorings. (Again, a coloring always means a proper vertex coloring.) Recall that the entropy argument in Proposition~\ref{prop:entropy_coloring} was not enough to give an asymptotically tight upper bound on $c_q(G)$. On the other hand, a combination of the entropy method and the proof idea for Lemma~\ref{lem:graph_containers} (rather than the statement itself) can produce sharp asymptotics for $c_q(G)$ for certain host graphs. For example, Kahn and the author obtained the asymptotics for the number of 4-colorings of $Q_n$. Recall that $N=2^n=|V(Q_n)|$. Throughout the section, we assume $n$ is large, and all the asymptotic notation is as $n \rightarrow \infty$.

\begin{theorem}[Kahn, P. \cite{Kahn2020Number}]\label{thm:KP_color}
\begin{equation}\label{eq:C4} c_4(Q_n)= 6e(1+o(1))2^{N}.\end{equation}
\end{theorem}

As is usual in the current article, the value on the right-hand side of \eqref{eq:C4} is an asymptotic lower bound on $c_4(Q_n)$: to see this, say the set of given colors is $\{1,2,3,4\}$. Consider a ``ground state" (\textit{cf.} Section~\ref{Subsubsec:lb_iQd}) that consists of the colorings of $V(Q_n)$ where all the even vertices take colors from $\{1,2\}$ and all the odd vertices take colors from $\{3,4\}$. Obviously, this construction produces $2^N$ distinct (proper) colorings. It is easy to verify that (similarly to the computation in Section~\ref{Subsubsec:lb_iQd}) the contribution of the colorings with non-nearby ``flaws" -- the vertices whose colors do not observe the ground state -- together with the colorings in the ground state is asymptotically $e\cdot 2^N$. Now the asymptotic lower bound follows by observing that there are six ground states.

So the main task for the proof of Theorem~\ref{thm:KP_color} is to show that the number of colorings with nearby (i.e., 2-linked) flaws is negligible (\textit{cf.} Section~\ref{sec:iQd_ub}). The problem of estimating $i(Q_n)$, the number of independent sets in $Q_n$, was able to be done by invoking Lemma~\ref{lem:graph_containers}. However, this is not the case for counting colorings. In the next example, we illustrate why the graph container method alone is not sufficient here.

\begin{example}\label{ex:cQd}
    Write $f$ for a member of $\mathcal C_q(Q_n)$ (that is, $f$ is a $q$-coloring of $V(Q_n)$) and let $[q]$ be the set of colors (palette). For simplicity, assume $q$ is even, and let $P=[1, q/2]$ and $Q=[q/2+1,q]$ (so $P$ and $Q$ equipartition the palette). Similarly to Example~\ref{ex:iQd}, we consider a special type of bad coloring. Let $\mathcal C \subseteq \mathcal C_q(Q_n)$ be the collection of $f$ such that
\begin{equation}\label{eq:colors} f(v) \in \begin{cases}
P \quad \text{ if $v \in \mathcal E$;}\\
Q \quad \text{ if $v \in \mathcal O$.}
\end{cases}
\end{equation}
(this will be our ground state), noticing that $|\mathcal C|=(q/2)^N$. Similarly to $\mathcal I_k$ in Example~\ref{ex:iQd}, define $\mathcal C_k \subseteq \mathcal C$ to be the collection of $q$-colorings such that the set $A$ of ``flaws" (meaning the vertices that violate \eqref{eq:colors}) in $\mathcal O$ is 2-linked and $|[A]|=k~( \in [n^4, N/4])$. For the sake of a cleaner exposition, let's further impose the (very strong) assumption that there are no flaws in $\mathcal E$.
\begin{figure}[h]
\begin{center}
    \includegraphics[height=2.5cm]{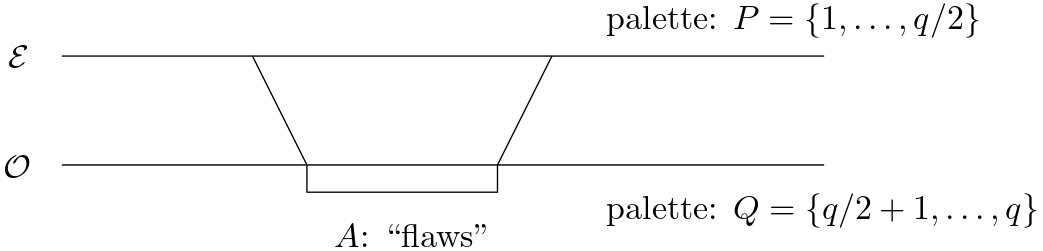}
\end{center}
\caption{Setting for Example~\ref{ex:cQd}}
\end{figure}

Lemma~\ref{lem:graph_containers} bounds the number of choices for $A$. Given $A \subseteq \mathcal O$, let $\mathcal C_A$ be the collection of $f \in \mathcal C_k$ whose set of flaws is precisely $A$. Observe that we have a bound
\[|\mathcal C_A|\le (q/2)^{N-|N(A)|}(q/2-1)^{|N(A)|}\]
(this is because if $w \in N(A)$, then $w$ is adjacent to a vertex $v \in A$ and the color of $v$ is in $P$, so $f(w)$ has at most $q/2-1$ choices (rather than $q/2$ choices)), and so
\[\begin{split}|\mathcal C_k|\le \sum_{v \in \mathcal E} \sum_{g}\sum_{\substack{A \subseteq \mathcal O, v \in N(A)\\ A \text{ 2-linked}, |[A]|=k, |N(A)|=g}} |\mathcal C_A|.\end{split}\]
The inner sum of the above is, by applying Lemma~\ref{lem:graph_containers} (and using the isoperimetric property of $Q_n$ that $|N(A)|-|A|=\Omega(|N(A)|/\sqrt n)$ as discussed in Example~\ref{ex:iQd}),
\[\begin{split}\sum_{\substack{A \subseteq \mathcal O, v \in N(A)\\ A \text{ 2-linked}, |[A]|=k, |N(A)|=g}} |\mathcal C_A|& \le  |\mathcal G(k, g, v)|\cdot (q/2)^{N} \left(\frac{q/2-1}{q/2}\right)^{g}\\
&\le (q/2)^N\cdot \left(\frac{2(q/2-1)}{q/2}\right)^g \cdot 2^{-\Omega(g/(\sqrt n\log n))}.\end{split}\]
Here the point is that the above expression is not necessarily smaller than the ``benchmark" $(q/2)^N$ if $q > 4$ (the case of $q=4$ is also easily seen to be problematic once we consider less specialized types of bad colorings). This issue didn't exist in enumerating independent sets, since, crucially, there a vertex has only two possible options for its ``color" (either occupied or unoccupied). However, as the number of options for the colors of vertices increases, there tends to be less ``penalty" from the specified flaws (in the above computation, the penalty from each vertex in $N(A)$ is $(q/2-1)/(q/2)=1-2/q$). So in such cases, the penalty is not enough to compensate the cost for flaws given by Lemma~\ref{lem:graph_containers}.
\end{example}

Kahn and the author \cite{Kahn2020Number} handled the above issue in proving Theorem~\ref{thm:KP_color} by carefully combining the entropy method and the graph container method. A somewhat vague description of their idea is as the following: Example~\ref{ex:cQd} shows that specifying flaws and identifying the colorings observing the given flaws in two separate steps is too expensive. On the other hand, just as independent sets in the given graph cluster together, colorings of the given graph also should ``cluster." As in this intuition, Kahn and the author regarded the containers in the proof (not in the statement) of Lemma~\ref{lem:graph_containers} as ``containers of colorings," and specified colorings directly from the containers. This way, they could skip the expensive step of specifying flaws. But of course, specifying colorings from containers requires much more involved analysis than specifying colorings given that the flaws are completely specified, and the entropy techniques were crucially used for this treatment.

Kahn and the author's approach for Theorem~\ref{thm:KP_color} was vastly developed by Jenssen and Keevash \cite{Jenssen2023Homomorphisms} and applied to much more general weighted graph homomorphisms from the discrete torus $\mathbb Z_m^n$, where $m$ is even, to any fixed graphs (and more). Their proof combines various techniques including entropy, graph containers, the cluster expansion method, and isoperimetric and algebraic properties of the torus.

Balogh, Garcia, and Li \cite{Balogh2021Independent} suggested the problem of enumerating colorings of the middle layers of the hypercube $Q_n$ for odd $n$. Note that when $n$ is odd, $Q_n$ has two maximum layers -- the $\frac{n\pm 1}{2}$-th layers. The subgraph of $Q_n$ induced by those two layers forms a $d$-regular graph where $d=(n+1)/2$. While this graph is still regular and satisfies some expansion properties (but the expansion is weaker than $Q_n$), its structure is not as ``well-organized" as $Q_n$ itself. In particular, some of the key arguments in \cite{Jenssen2023Homomorphisms} based on the symmetric structure of $Q_n$ does not immediately extend to the middle layers of $Q_n$. Li, McKinley, and the author \cite{Li2025Number} addressed the question of Balogh, Garcia, and Li when the number of colors is even. They improved Kahn and the author's approach so that, as long as the number of colors is even, the combination of the entropy method and the container approach works without taking advantage of the nice symmetries in the structure of the hypercube. They noted that their proof can be adapted for a broader class of graphs beyond the middle layers of $Q_n$. The case where the number of colors is odd remains open.

\subsection{$\mathbb Z$-homomorphisms from the hypercube}\label{sec:Z-hom}

Probably somewhat surprisingly, the problem of enumerating colorings of the hypercube was motivated by the study of the typical ``range" of $\mathbb Z$-homomorphisms from $Q_n$. The solutions of this problem again involve a combination of entropy and graph containers, but here the use of the entropy method is quite different from that for the proof of Theorem~\ref{thm:KP_color}.

We start with some relevant definitions. For a finite, connected, bipartite graph $G=(V,E)$ with distinguished vertex $v_0 \in V$, set
\[\Hom_{v_0}(G)=\{f:V \rightarrow \mathbb Z: f(v_0)=0, \{x,y\} \in E \Rightarrow |f(x)-f(y)|=1\}.\]
That is, $\Hom_{v_0}(G)$ is the collection of $\mathbb Z$-homomorphisms $f$ from $G$ with $f(v_0)=0.$ Note that $\Hom_{v_0}(G)$ is a finite set. For $f \in \Hom_{v_0}(G)$, define the \textit{range} of $f$ by
\[R(f):=|\{f(x):x \in V\}|.\]
For any set $X$, we write $x \in_r X$ to denote that $x$ is a uniformly randomly chosen element of $X$.

Benjamini, H\"aggstr\"om, and Mossel \cite{Benjamini2000Random} investigated the distribution of $R(f)$ for $f \in_r \Hom_{v_0}(G)$ for various host graphs $G$, with a motivating question of whether concentration inequalities for \textit{typical} Lipschitz functions are stronger than those which hold for \textit{all} Lipschitz functions. We refer the readers to \cite{Benjamini2000Random} for more details of this interesting direction of research, and here we just focus on the hypercube. 

Observe that, when the host graph is $Q_n$, $R(f)$ can be as large as $n+1$ (we may assume $f(\bar 0)=0$, and then assign $f \equiv i$ to the $i$th layer for $i \in [1, n]$). However, for $f \in_r \Hom_{v_0}(Q_n)$, it should be highly unlikely that $f$ achieves such a large range. Indeed, it was conjectured in \cite{Benjamini2000Random} that if $f\in_r \Hom_{v_0}(Q_n)$, then $R(f)=o(n)$ w.h.p. This conjecture was proved in a very strong way by Kahn~\cite{Kahn2001Range} using a beautiful entropy argument.

\begin{theorem}[Kahn \cite{Kahn2001Range}] \label{thm:Kahn_cube}
There is a constant $b$ for which the following holds. For $f \in_r\Hom_{v_0}(Q_n)$,
\[\mathbb P(R(f)>b)<2^{-\Omega(n)}.\]    
\end{theorem}
In other words, $R(f)=O(1)$ w.h.p. Kahn further conjectured that, in fact, the constant should be $5$, and this conjecture was proved by Galvin~\cite{Galvin2003Homomorphisms}. To see why $5$ should be the answer: note that a ``ground state" for $\mathbb Z$-homomorphisms from $Q_n$ is to assign 0 to one part (say, $\mathcal E$) and let each of the odd vertices choose a value from $\{\pm 1\}$. This gives $R(f)=3$. However, similarly to the discussion in Section~\ref{Subsubsec:lb_iQd}, it is easy to see that planting some non-nearby flaws in $\mathcal E$ and assigning them the values $\pm 2$ produces nontrivial contributions to the homomorphism count. In other words, having $R(f)=5$ is also not unlikely (in fact, the probabiliy for this event is $\Omega(1)$).
\begin{figure}[h]
\begin{center}
\includegraphics[height=3.5cm]{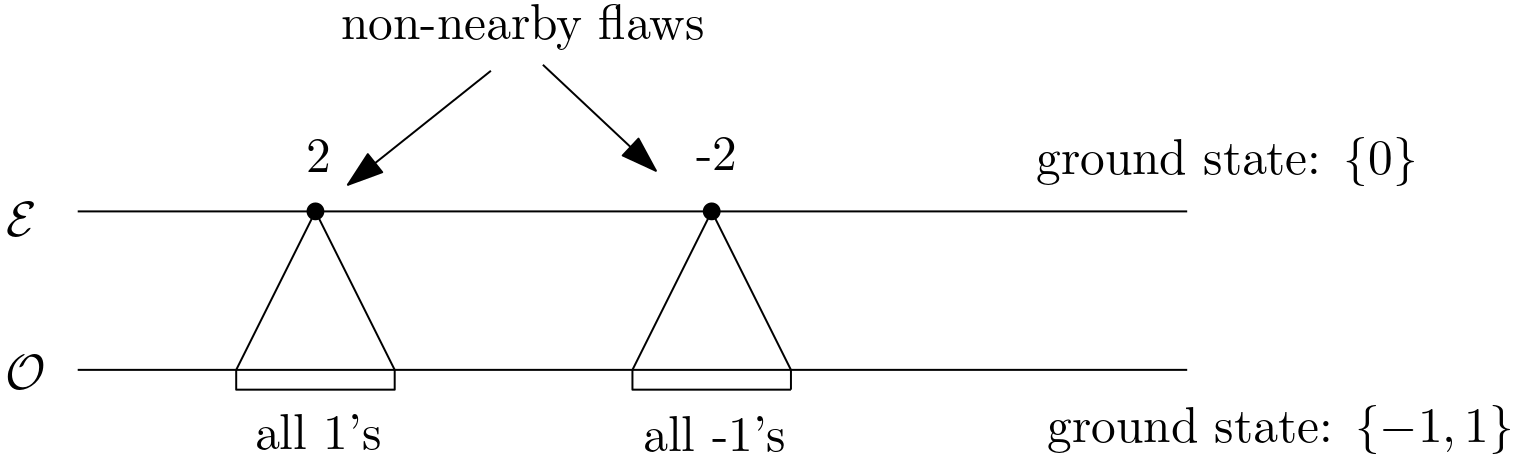}
\end{center}
\caption{A configuration of $f$ with $R(f)=5$}
\end{figure}

\noindent On the other hand, having a vertex with the value greater than or equal to $\pm 3$ requires a large set of nearby flaws, which should be highly unlikely.

The starting point of Galvin's proof is Theorem~\ref{thm:Kahn_cube}, and he used the graph container method to find sharp asymptotics for $|\Hom_{v_0}(Q_n)|$ and $|\{f \in \Hom_{v_0}(Q_n):R(f)\le 5\}|$. By comparing the two asymptotics, the conclusion that the $b$ in Theorem~\ref{thm:Kahn_cube} should be 5 follows as an immediate corollary. Peled~\cite{Peled2017High} vastly generalized Galvin's result, obtaining results for a general class of tori including $\mathbb Z_n^d=(\mathbb Z/n\mathbb Z)^d$ (where $n$ can be large with respect to $d$) and providing strong bounds on $\mathbb P(R(f) \ge k)$ for arbitrary $k$, for both random $\mathbb Z$-homomorphisms and 1-Lipschitz functions.

D. Randall made the interesting observation that there is a bijection from the set $\Hom_{v_0}(Q_n)$ to  the set of proper $3$-colorings of $V(Q_n)$ with the color of $v_0$ fixed. The bijection is given by $f \rightarrow \chi$ where $\chi(v)=i$ iff $f(v) \equiv i$ (mod 3). As a result, Galvin's enumeration of $\mathbb Z$-homomorphisms provides the sharp asymptotics for the number of $3$-colorings of $Q_n$. Later, Engbers and Galvin \cite{Engbers2012Hcoloring} extended Kahn's entropy approach for Theorem~\ref{thm:Kahn_cube} to general graph homomorphisms from the discrete torus $\mathbb Z_m^n$, where $m$ is even, to any fixed graph. This work was a starting point of Theorem~\ref{thm:KP_color} and inspired the aforementioned work of Jenssen and Keevash \cite{Jenssen2023Homomorphisms}.

\subsection{Lipschitz functions from expanders}\label{sec:Lips}

Motivated by the study of the range of random $\mathbb Z$-homomorphisms on tree-like graphs (\cite{Benjamini2000Random}), hypercubes (\cite{Kahn2001Range, Galvin2003Homomorphisms}), and $\mathbb Z^d$ for large $d$ (\cite{Peled2017High}), the typical range of Lipschitz functions on expanders was investigated by Peled, Samotij, and Yehudayoff~\cite{Peled2013Lipschitz}. We first quickly introduce some relevant definitions.

Recall that, for given $G=(V,E)$, an integer-valued function $f:V(G)  \rightarrow \mathbb Z$ is called an \textit{$M$-Lipschitz function} from $G$ if $|f(u)-f(v)| \le M$ for all $\{u,v\} \in E(G)$. For a vertex $v_0 \in V(G)$, let
\[\Lip_{v_0}(G;M)=\{f:V(G) \rightarrow \mathbb Z~|~ \text{$f(v_0)=0$ and $|f(u)-f(v)|\le M$ $ \forall \{u,v\} \in E(G)$}\}.\]
Define the \textit{range} of $f \in \Lip_{v_0}(G;M)$ (\textit{cf.} Section~\ref{sec:Z-hom}) by
\[R(f):=\max_{v \in V(G)}f(v)-\min_{v \in V(G)}f(v)+1.\]

The main focus of this section is the typical range of Lipschitz functions defined on expander graphs. Following \cite{Peled2013Lipschitz}, say a $d$-regular $N$-vertex graph $G$ is a \textit{$\lambda$-expander} if it satisfies
\[\text{$\left|e(S,T)-\frac{d}{N}|S||T|\right|\le \lambda \sqrt{|S||T|}$ for all $S, T \subseteq  V(G)$.}\]
(We always have $(1-d/n)\sqrt d \le \lambda \le d$ (see, e.g., \cite[Section 1.1]{Peled2013Lipschitz}), and the discussion in this section focuses on $\lambda \le d/5$.) We remark that the above definition of $\lambda$-expander (which is motivated by the {expander mixing lemma} \cite{Alon1988Explicit}) guarantees some nice expansion properties for $G$ as its name indicates: for example, for any $v \in V(G)$ and $t \in \mathbb Z^+$, the volume of the ball in $G$ of radius $t$ centered at $v$ (i.e., the set of vertices $w$ with $\dist_G(v,w) \le t$) is at least $\min\left\{\frac{N}{2}, \left(\frac{d}{2\lambda}\right)^{2t}\right\}$. So roughly it says that as long as $t$ is ``not too large," the volume of the ball grows exponentially. This volume growth poses an upper bound on the range of $M$-Lipschitz functions from $G$: it can be derived from the volume growth that any $\lambda$-expander $G$ has diameter ${\rm diam}(G) \le \log\left(\frac{d}{2\lambda}\right)^{-1}\log N$, thus trivially, any $f \in \Lip_{v_0}(G;M)$ must satisfy
\begin{equation}\label{eq:range} R(f)\le M\cdot\log\left(\frac{d}{2\lambda}\right)^{-1}\log N+1.\end{equation}
On the other hand, it was shown in \cite{Peled2013Lipschitz} that if $G$ is a ``sufficiently good" expander (roughly, the authors of \cite{Peled2013Lipschitz} require 
\begin{equation}\label{eq:expansion}|N(A)| \ge \Omega( M\log(dM)|A|)\end{equation}
for each ``not-too-large" set $A \subseteq V(G)$) then for $f \in_r \Lip_{v_0}(G;M)$, w.h.p. $R(f)$ is much smaller than the trivial bound in \eqref{eq:range}; more precisely, they showed that the $\log N$ in \eqref{eq:range} can be replaced by $\log\log N$.

Consideration of $M$-Lipschitz functions for $M \rightarrow \infty$ is a natural bridge to \emph{$\mathbb R$-valued Lipschitz functions}, functions $f: V(G) \to \mathbb{R}$ satisfying $|f(u)-f(v)|\le 1$ for all $\{u,v\} \in E(G)$. It is not hard to see that (e.g. \cite[page 8]{Peled2013Grounded}) if $f \in_r \Lip_{v_0}(G;M)$, then as $M \rightarrow \infty$, $f/M$ converges in distribution to a uniformly random $\mathbb R$-valued Lipschitz function on $G$ with $f(v_0)=0$. With this motivation, it is natural to ask whether one can remove the dependence on $M$ in the expansion requirement \eqref{eq:expansion} and still obtain a similar conclusion as in \cite{Peled2013Lipschitz}. Indeed, it was shown by Krueger, Li, and the author \cite{Krueger2024Lipschitz} that requiring an $\Omega(1)$-expansion (i.e., $|N(A)| \ge \Omega(|A|)$) is enough as long as $M$ satisfies the condition in the below statemet:
\begin{theorem}[Krueger-Li-P. \cite{Krueger2024Lipschitz}]\label{cor.flat}
There exist universal constants $c, C, c', C'>0$ such that the following holds. Let $M \in \mathbb Z^+$, and let $G$ be a connected $N$-vertex, $d$-regular $\lambda$-expander with $\lambda \le d/5$, where $M \le \min\{c\frac{d^{3/2}}{\lambda \log d}, (\log N)^C\}$. Let  $f$ be a uniformly randomly chosen $M$-Lipschitz function on $G$. Then
\begin{equation}\label{eq:Lips} \mathbb P\left( R(f) \geq C'M \frac{\log\log N}{\log(d/\lambda)} +2(M+1) \right) \leq N^{-c'}. \end{equation}
\end{theorem}

\noindent (to read off the $\Omega(1)$-expansion from the theorem above, we note that the expansion for $\lambda$-expanders is $\Omega(d/\lambda)$). 

It was noted in~\cite{Peled2013Lipschitz} that the lower bound on $R(f)$ in \eqref{eq:Lips} should be tight, and one may apply the techniques of Benjamini, Yadin, and Yehudayoff \cite{Benjamini2007Random} to prove it. This idea was implemented (with some additional ideas for the case $M \gg 1$) by I\c{s}\i k and the author \cite{Isik2024Random}. 

The proof of Theorem~\ref{cor.flat} is motivated by the idea of Li, McKinley, and the author \cite{Li2025Number} in the study of random colorings of the middle layers of the hypercube. The restriction on $M$ in Theorem~\ref{cor.flat} was due to the limit of the techniques used in \cite{Krueger2024Lipschitz}. As mentioned above, the study of $M$-Lipschitz functions for $M \rightarrow \infty$ provides a transition to $\mathbb{R}$-valued Lipschitz functions, and a natural question is whether $\mathbb{R}$-valued Lipschitz functions still exhibit a small range (much smaller than the diameter of the host graph) on expander graphs (and in what quantitative sense). This question is still wide open, and a possible first step towards this direction seems to be to understand the behavior of $M$-Lipschitz functions for $M \gg d$.

\section{Antichains in $[t]^n$}\label{sec:t^n}

This section discusses a first successful application of the graph container lemma for expanders in the ``irregular" setting. Recall that Dedekind's problem concerns the number of antichains in the Boolean lattice. A natural generalization of the Boolean lattice is the poset $([t]^n, \preceq)$ (we will simply write $[t]^n$) -- the ground set consists of all $n$-tuples $(x_1,\dots,x_n)$ of integers in $\{0,1,\dots,t-1\}$, and the partial order $\preceq$ is defined to be $x\preceq y\Leftrightarrow x_i\leq y_i$ for all $i \in [n]$. We write $\alpha([t]^n)$ for the number of antichains in the poset $[t]^n$.

\subsection{Connection to Ramsey-type problems and logarithmic asymptotics}\label{sec.Ramsey} Enumerating antichains in $[t]^n$ is a natural generalization of Dedekind's problem, but it also turns out to have an interesting connection to a Ramsey type question. This connection is not directly related to the focus of the present article so we will try to be brief and refer the readers to \cite{FalgasRavry2023Dedekinds, Fox2012Erdos, Moshkovitz2012Ramsey} for details and omitted definitions. Fox, Pach, Sudakov, and Suk \cite{Fox2012Erdos} introduced the number $M_k(t,n)$ as the smallest integer $M$ with the property that for any $n$-colorings of ${[M] \choose k}$ one can find a monochromatic ``monotone path"\footnote{Formally, for any sequence of positive integers $x_1<x_2<\ldots<{x_{l+k-1}}$, we say that the $k$-tuples $(x_i,x_{i+1}, \ldots, x_{i+k-1})$ $(i=1,2,\ldots, l)$ form a \textit{monotone path} of length $l$.} of length $t$. The celebrated results of Erd\H{o}s and Szekeres \cite{Erdos1935Combinatorial} can be expressed using $M_2(t,n)$ and $M_3(t,2)$, and Fox et al suggested the problem of estimating $M_k(t,n)$ for general $k, t$ and $n$. Moshkovitz and Shapira \cite{Moshkovitz2012Ramsey} and Milans, Stolee, and West \cite{Milans2015Ordered} independently observed that
\[\alpha([t]^n)=M_3(t,n)-1 \quad \text{for all }  t,n \ge 2,\]
and using this relationship, Moshkovitz and Shapira obtained the following bounds on $M_3(t,n)$:
\begin{theorem}[Moshkovitz, Shapira \cite{Moshkovitz2012Ramsey}] For every $t,n \ge 2$,
\[2^{\frac{2}{3}t^{n-1}/\sqrt n} \le M_3(t,n) \le 2^{2t^{n-1}}.\]
\end{theorem}
Moshkovitz and Shapira boldly suggested that, as $n \rightarrow \infty$,
\begin{equation}\label{MS.bold} M_3(t,n)~ ( =\alpha([t]^n)+1 )~=2^{(1+o_n(1))N(t,n)}\end{equation}
where $N(t,n)$ is the size of a largest layer of $[t]^n$. Notice the resemblance of the formulation in \eqref{MS.bold} to Theorem~\ref{thm:Kleitman}. After partial progress towards \eqref{MS.bold} with some restrictions on the ranges of $t$ and $n$ (e.g., \cite{Carroll2012Counting, Noel2018Supersaturation, Tsai2019Simple, Pohoata2021Number, Park2024Note}), Falgas-Ravry, R\"aty and Tomon \cite{FalgasRavry2023Dedekinds} finally confirmed \eqref{MS.bold} in full, using the container method based on a cleverly constructed normalized matching flow:
\begin{theorem}[Falgas-Ravry, R\"aty and Tomon \cite{FalgasRavry2023Dedekinds}]\label{thm:Falgas}
There exists a constant $c>0$ such that for every $n, t \ge 2$, we have
\begin{equation}\label{eq:log.tn} \log_2\alpha([t]^n)\le \left(1+c\cdot \frac{(\log n)^3}{n}\right)\cdot N(t,n).\end{equation}
\end{theorem}

\subsection{Towards sharp asymptotics}

Concerning the sharpness of the error term in \eqref{eq:log.tn}, Falgas-Ravry et al \cite{FalgasRavry2023Dedekinds} showed that
\begin{equation}\label{eq:FRT_conj}\log_2\alpha([t]^n)\ge \left(1+2^{-\Theta(n)}\right)\cdot N(t,n),\end{equation}
and said that ``it might be reasonable to conjecture that this lower bound is closer to optimal than our upper bound (in Theorem~\ref{thm:Falgas})." For $t=2$, an upper bound matching the right-hand side of \eqref{eq:FRT_conj} follows from Theorem~\ref{thm:Dedekind}, so one might expect that the proof techniques for Dedekind's problem might extend to $\alpha([t]^n)$, at least for $t=O(1)$. While Lemma~\ref{lem:graph_containers} works well for Dedekind's problem, somewhat surprisingly, this lemma is not immediately usable for counting antichains in $[t]^n$ for $t \ge 3$. In fact, Sapozhenko \cite{Sapozhenko2005Systems} himself asked what the number of antichains in $[3]^n$ is. The question of estimating $\alpha([3]^n)$ was also asked by Noel, Scott and Sudakov \cite{Noel2018Supersaturation}, who proved asymptotics for $\log_2\alpha([3]^n)$.

We briefly explain why Sapozhenko's approach for Dedekind's problem doesn't immediately extend to $\alpha([3]^n)$. Recall that, in the graph setting, Sapozhenko's argument crucially assumes the degree condition in \eqref{eq:Sap_deg}:
    \[ \max_{v \in Y} \deg(v) \le \min_{v \in X} \deg(v).\]
Note that the graphs induced by any two consecutive layers of the Boolean lattice $\mathcal B_n$ satisfy the above degree condition (by choosing the parts $X$ and $Y$ appropriately), and this was a crucial ingredient for Sapozhenko's (asymptotic) resolution of Dedekind's problem. On the other hand, the graphs induced by two consecutive layers of $[t]^n$ do not necessarily satisfy the degree condition in \eqref{eq:Sap_deg} when $t \ge 3$.

Sapozhenko \cite{Sapozhenko2005Systems} mentioned that the irregular structure of $[3]^n$ is the main obstacle for the application of his method. Recently, Jenssen, Sarantis, and the author \cite{Jenssen2025Number} found a way to handle this obstacle and obtained the following sharp estimate for $\alpha([3]^n)$. Their approach is still based on the graph container method, but it is the first to extend Lemma~\ref{lem:graph_containers} beyond the degree requirement in \eqref{eq:Sap_deg}.

\begin{theorem}[Jenssen, P., Sarantis \cite{Jenssen2025Number}]\label{thm:3n} Let $N(3,n)$ be the size of the middle layer of $[3]^n$. Then,
$$\alpha([3]^n)=2^{N(3,n)}\exp\left[(1+o(1))\sqrt{\frac{1+2\sqrt{2}}{2\sqrt{2}\pi}}n^{-1/2}\left(\frac{1+2\sqrt{2}}{2}\right)^n\right].$$
\end{theorem}

As we have discussed several times in the previous sections, the main result of \cite{Jenssen2025Number} is structural, and the quantitative bound in Theorem~\ref{thm:3n} follows as a corollary from the structural result. The statement about the structure of almost all antichains in $[3]^n$ is  very similar to Theorem~\ref{thm:Dedekind_structure}, and the estimate in Theorem~\ref{thm:3n} follows by analyzing the contributions of antichains with small flaws.

We note that, for $t \ge 4$, the problem of obtaining a similarly sharp estimate for the number of antichains in $[t]^n$ is still wide open, and we don't even know whether the lower bound in \eqref{eq:FRT_conj} is optimal. It might be natural to suspect that for bounded $t$, a structural result similar to Theorem~\ref{thm:Dedekind_structure} still holds, which would give a much stronger answer to the question of Falgas-Ravry et al following \eqref{eq:FRT_conj}.

\section{Maximal independent sets}\label{sec:mis}

Sometimes applications of the graph container method lead to an interesting development of isoperimetric inequalities on the host graphs. One nice example is the problem of counting maximal independent sets in the hypercube.

An independent set $I \in \mathcal I(G)$ is called \textit{maximal} if for any $v \in V(G)$, $I \cup \{v\}$ is no longer independent. Write ${\rm MIS}(G)$ for the family of maximal independent sets in $G$, and $\mis(G):=|{\rm MIS}(G)|$. In general, finding $\mis(G)$ is $\texttt{\#} P$-hard, but its extremal properties are well-known: for example, Moon and Moser \cite{Moon1965Cliques} showed that $\mis(G) \le 3^{N/3}$ for any $N$-vertex graph $G$, and the unique extremal graph is a vertex disjoint union of $N/3$ triangles if $3$ divides $N$. Hujter and Tuza \cite{Hujter1993Number} showed that if $G$ is triangle-free, then the upper bound is reduced to $\mis(G) \le 2^{N/2}$, and the unique extremal graph is a matching with $N/2$ edges if $2$ divides $N$.

The problem of enumerating maximal independent sets in $Q_n$ was asked by Duffus, Frankl, and R\"odl \cite{Duffus2013Maximal}. Observe that (recalling that $N=|V(Q_n)|=2^n$)
\[\mis(Q_n) \le 2^{N/2};\]
in fact, for any bipartite graph $G$ with a bipartition $X \cup Y$, $\mis(G) \le \max\{2^{|X|}, 2^{|Y|}\}$, since if $I \in {\rm MIS}(G)$, then $I \cap Y=Y \setminus N(I \cap X)$ (so $I \cap X$ determines $I \cap Y$). On the other hand, by considering a certain construction of maximal independent sets (\textit{cf.} Section~\ref{Subsubsec:lb_iQd}), one can obtain the following asymptotic lower bound.
\begin{observation}\label{obs:misQd}
    \[{\rm mis}(Q_n)>(1-o(1))2n2^{N/4}.\]
\end{observation}

We briefly sketch a justification of the above observation. (For a full justification, see \cite[Section 1.1]{Kahn2022Number}.) An \textit{induced matching} (IM) in a graph is an induced subgraph that is a matching, and a maximum IM in $Q_n$ is constructed as the following. Pick an $i \in [n]$, and consider $Q_0:=\{v \in V(Q_n): v_i=0\}$ and $Q_1=\{v \in V(Q_n): v_i=1\}$ where $v_i$ denotes the $i$-th coordinate of $v$.
\begin{figure}[h]
\begin{center}
\includegraphics[height=6cm]{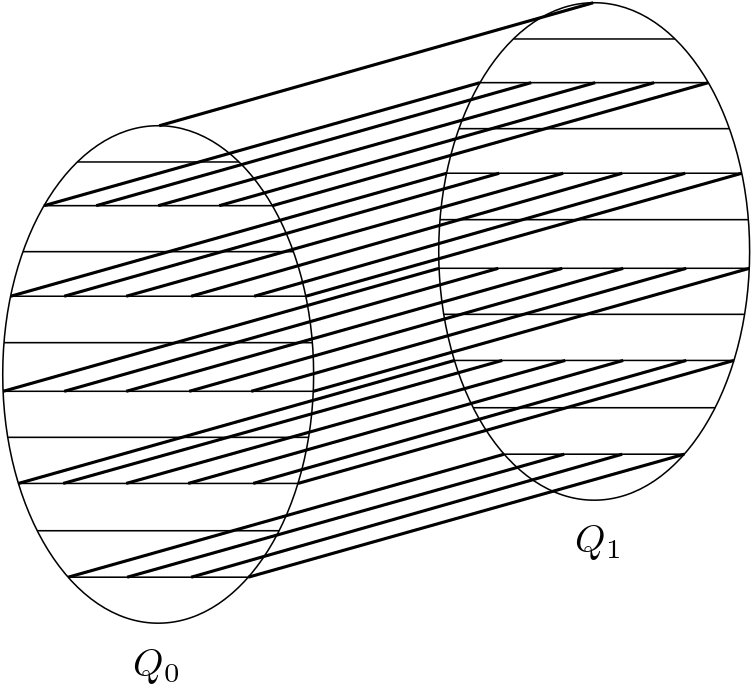}
\caption{A maximum induced matching between $Q_0$ and $Q_1$}
\end{center}
\end{figure}

Note that $\nabla(Q_0, Q_1)$ is a (perfect) matching. Now, pick the edges of $\nabla(Q_0, Q_1)$ with its end in $V(Q_0)$ having an odd Hamming weight, then these edges form an IM of $Q_n$ with $N/4$ edges. It is easy to see that they are maximum IM's in $Q_n$, and there are $2n$ of such IM's ($n$ choices for $i$ and 2 choices for the parity of the ends in $V(Q_0)$). We leave it to the readers that those are the only maximum IM's.

Denote by $M^*$ the maximum IM's constructed above. Each $M^*$, which plays the role of a ``ground state" in the MIS count, gives rise to exactly $2^{N/4}$ MIS's, gotten by choosing one vertex from each edge of $M^*$ and extending the resulting independent set to the (unique) MIS containing it. It is also easy to see that the number of MIS's that belong to more than one ground state are negligible. Since there are $2n$ ground states, the asymptotic lower bound in Observation~\ref{obs:misQd} follows.

\smallskip

It was conjectured by Ilinca and Kahn \cite{Ilinca2012Counting}, who proved that $\log\mis(Q_n) \sim N/4$ as conjectured by Duffus, Frankl, and R\"odl \cite{Duffus2013Maximal}, that $\mis(Q_n)$ is asymptotically equal to the lower bound in Observation~\ref{obs:misQd}. This conjecture, which again is from the framework in Section~\ref{sec:framework}, was proved by Kahn and the author \cite{Kahn2022Number}.

\begin{theorem}[Kahn, P. \cite{Kahn2022Number}]
    \[{\rm mis}(Q_n) \sim 2n2^{N/4}.\]
\end{theorem}

Along with the graph container method, a crucial ingredient for the above result is the isoperimetric inequality in Theorem~\ref{thm:isoper} below. The inequality was used to show that (we skip the derivation and refer to \cite{Kahn2022Number} for details) almost all MIS's have a large intersection with some maximum induced matching $M^*$. This step establishes the fact that each $M^*$ indeed forms a ground state, and so it suffices to analyze MIS's with small flaws in each ground state. 

In the following statement, $\mu$ is the uniform measure on $V=V(Q_n)$. For $x \in V$ and $T \subseteq V$, let $d_T(x)$ be the number of neighbors of $x$ in $T$ and define $h_S:V \rightarrow \mathbb N$ by
\[h_S(x)=
\begin{cases}
d_{V \setminus S}(x) & \mbox{ if $x \in S$,}\\
0 & \mbox{ if $x \notin S$.}
\end{cases}\]
We define $\beta=\log_2(3/2) (\approx .585)$ and $\int f d\mu = \sum_{x \in V} f(x) \mu(x).$ 

\begin{theorem}\label{thm:isoper}
For any $A \subseteq V$,
\[\int h_A^\beta d\mu \ge 2 \mu(A) (1-\mu(A)).\]
\end{theorem}
 We refer to \cite{Kahn2020Isoperimetric} for detailed discussion on the consequences of the above theorem, and here we just state a simple implication to hint at the flavor of Theorem~\ref{thm:isoper}. Define  $\nabla A:=\{\{a,b\} \in E(G):a \in A, b \notin A\}$.  

 \begin{cor}
     For $k=1,2$, if $|A|=2^{n-k}$ and $|\nabla A|<(1+\epsilon)|A|\log(2^n/|A|)$, then there is a subcube $C$ with $\mu(C \Delta A)=O(\epsilon)$ (where $C \Delta A$ denotes the symmetric difference between $C$ and $A$).
 \end{cor}

Another question considered in \cite{Duffus2013Maximal} was enumerating maximal independent sets in the graph induced by the middle layers of $Q_n$ for odd $n$, which we denote by $\mathcal M_n$. Again, there is an easy construction that shows that $\log\mis(\mathcal M_n) \ge {n-1 \choose n/2}$, and Ilinca and Kahn showed that this lower bound is asymptotically tight for $\log\mis(\mathcal M_n)$, proving a conjecture of Duffus, Frankl, and R\"odl \cite{Duffus2013Maximal}. Recently, Balogh, Chen and Garcia \cite{Balogh2025Maximal} obtained sharp asymptotics for $\mis(\mathcal M_n)$. The key ingredients of their proof again include the graph container method and isoperimetric inequalities for $\mathcal M_n$.



\thankyou{The author is grateful to the referee for their careful reading of the manuscript and their constructive comments. We also thank Jozsef Balogh and Matthew Jenssen for their helpful feedback on the draft. The author was supported by NSF Grant DMS-2324978, NSF CAREER Grant DMS-2443706 and a Sloan Fellowship.}



\myaddress


\end{document}